\documentclass[11pt]{article}

\newcommand{\documentdate}{8 X 2021}

\usepackage{a4wide,graphicx,latexsym,varioref,amsmath}

\newcommand{\al}[1]{{\footnotesize {\sf #1}}}
\newcommand{\eo}[1]{\calO\left(#1\right)}
\newcommand{\s}[1]{^{\mbox{\protect\tiny #1}}}
\newcommand{\kap}[1]{\kappa_{\mbox{\tiny #1}}}
\newcommand{\evbndf}[1]{\mbox{\textrm{\textup{\small\texttt{\#}}}}_{\mbox{\tiny {\sf #1}}}\s{F}}
\newcommand{\evbndd}[1]{\mbox{\textrm{\textup{\small\texttt{\#}}}}_{\mbox{\tiny {\sf #1}}}\s{D}}

\newcommand{\kat}[1]{\kap{\sf #1}}
\newcommand{\kata}[1]{\kat{#1}\s{A}}
\newcommand{\katb}[1]{\kat{#1}\s{B}}
\newcommand{\katc}[1]{\kat{#1}\s{C}}
\newcommand{\katd}[1]{\kat{#1}\s{D}}
\newcommand{\kate}[1]{\kat{#1}\s{E}}
\newcommand{\katf}[1]{\kat{#1}\s{F}}
\newcommand{\kats}[1]{\kat{#1}\s{S}}
\newcommand{\epsmin}{\epsilon_{\min}}
\newcommand{\sphi}{\widehat{\phi}}
\newcommand{\barphi}{\overline{\phi}}
\newcommand{\accuracy}{{\tt accuracy}}
\newcommand{\accuracys}{{\tt accuracy}$_s$\,}
\newcommand{\accuracyj}{{\tt accuracy}$_j$\,}

\newcommand{\sufficient}{{\tt sufficient}}
\newcommand{\insufficient}{{\tt insufficient}}
\newcommand{\absolute}{{\tt absolute}}
\newcommand{\relative}{{\tt relative}}
\newcommand{\barDT}{\overline{\Delta T}\!}
\newcommand{\ngap}{\!\!\!}

\newcommand{\numsection}[1]{\section{#1}\setcounter{equation}{0}}

\newcommand{\beqn}[1]{\begin{equation}\label{#1}}
\newcommand{\eeqn}{\end{equation}}
\newcommand{\req}[1]{(\ref{#1})}
\newcommand{\bpr}{{\bf Proof.} \hspace{1.5mm}}
\newcommand{\epr}{\hfill $\Box$ \vspace*{1em}}
\newcommand{\proof}[1]{
\begin{list}{}{
\setlength{\topsep}{0.0pt}
\setlength{\partopsep}{0.0pt}
\setlength{\leftmargin}{0.025\textwidth}
\setlength{\rightmargin}{0.5\leftmargin}
\setlength{\labelwidth}{0.5\leftmargin}
\setlength{\labelsep}{0.25\leftmargin}}
\item \bpr #1 \epr \noindent
\end{list}}
\newcommand{\ms}{\;\;\;\;}
\newcommand{\tim}[1]{\;\; \mbox{#1} \;\;}
\newcommand{\bctable}[1]{\begin{table}[htbp]
                         \begin{center}
                         \begin{tabular}{#1} }
\newcommand{\ectable}[1]{\end{tabular}
                         \caption{#1}
                         \end{center}
                         \end{table}}
\newtheorem{theorem}{Theorem}[section]
\newtheorem{lemma}[theorem]{Lemma}
\newtheorem{corollary}[theorem]{Corollary}

\newcommand{\llem}[2]{\vspace{\baselineskip} 
\noindent\framebox[\textwidth]{\parbox{0.95\textwidth}{
\begin{lemma} \label{#1} \rm #2 \end{lemma} } } \vspace{\baselineskip} }
\newcommand{\lthm}[2]{\vspace{\baselineskip} 
\noindent\framebox[\textwidth]{\parbox{0.95\textwidth}{
\begin{theorem} \label{#1} \rm #2 \end{theorem} } } \vspace{\baselineskip} }

\newlength{\thmw}
\newcommand{\lbthm}[3]{\vspace{\baselineskip}\noindent\hbox{%
  \lower\fboxrule\hbox{\vbox{\hrule\hbox{\vrule \kern-\fboxrule \vbox{%
  \vspace{\fboxsep} \noindent\hspace{2\fboxsep}\parbox{\thmw}{
  \begin{theorem}\label{#1}{\rm #2}\end{theorem}\vspace{-\lastskip}}
  \hspace{\fboxsep}}\kern-\fboxrule \vrule }}}}\newpage \hbox{%
  \lower\fboxrule\hbox{\vbox{\hbox{\vrule \kern-\fboxrule \vbox{%
  \noindent\hspace{2\fboxsep}\parbox{\thmw}{\rm #3}\hspace{\fboxsep}
  \vspace{4\fboxsep}}\kern-\fboxrule \vrule }\hrule }}}\vspace{\baselineskip}
}
 
\newcommand{\calE}{{\cal E}}
\newcommand{\calO}{{\cal O}} 
\newcommand{\calS}{{\cal S}} \newcommand{\calU}{{\cal U}}
\newcommand{\sfrac}[2]{{\scriptstyle \frac{#1}{#2}}}
\newcommand{\half}{\sfrac{1}{2}}
\newcommand{\quarter}{\sfrac{1}{4}}
\newcommand{\comment}[1]{}

\newcommand{\eqdef}{\stackrel{\rm def}{=}}
\newcommand{\bigmin}{\displaystyle \min}
\newcommand{\bigmax}{\displaystyle \max}
\newcommand{\bigfrac}[2]{\frac{\displaystyle #1}{\displaystyle #2}}
\newcommand{\bigsum}{\displaystyle \sum}
\newcommand{\ii}[1]{\{1, \ldots, #1 \}}
\newcommand{\iiz}[1]{\{0, \ldots, #1 \}}

\renewcommand{\Re}{\hbox{I\hskip -2pt R}}

\newcounter{algo}[section]
\renewcommand{\thealgo}{\thesection.\arabic{algo}}
\newcommand{\algo}[3]{\refstepcounter{algo}
\begin{center}\begin{figure}[htbp]
\framebox[\textwidth]{
\parbox{0.95\textwidth} {\vspace{\topsep}
{\bf Algorithm \thealgo : #2}\label{#1}\\
\vspace*{-\topsep} \mbox{ }\\
{#3} \vspace{\topsep} }}
\end{figure}\end{center}}
\newcommand{\barf}{\overline{f}}
\newcommand{\barT}{\overline{T}}

\title{Strong Evaluation Complexity of An Inexact Trust-Region Algorithm with
       for Arbitrary-Order Unconstrained Nonconvex Optimization}

\author{
C. Cartis\thanks{Mathematical Institute,
   Oxford University,
   Oxford OX2 6GG, England.  Email: coralia.cartis@maths.ox.ac.uk},
N. I. M. Gould\thanks{Computational Mathematics Group,
   STFC-Rutherford Appleton Laboratory,
   Chilton OX11 0QX, England. Email:  nick.gould@stfc.ac.uk .
   The work of this author was supported by EPSRC grant EP/M025179/1}
~and~Ph. L. Toint\thanks{Namur Center for Complex Systems (naXys),
   University of Namur, 61, rue de Bruxelles, B-5000 Namur, Belgium.
   Email: philippe.toint@unamur.be}
}

\date{\documentdate}

\begin{document}

%  The manuscript

\maketitle

\begin{abstract}
A trust-region algorithm using inexact function and derivatives values is
introduced for solving unconstrained smooth optimization problems. This
algorithm uses high-order Taylor models and allows the search of strong
approximate minimizers of arbitrary order.  The evaluation complexity of
finding a $q$-th approximate minimizer using this algorithm is then shown,
under standard conditions, to be
$\mathcal{O}\big(\min_{j\in\ii{q}}\epsilon_j^{-(q+1)}\big)$ where the
$\epsilon_j$ are the order-dependent requested accuracy
thresholds. Remarkably, this order is identical to that of classical
trust-region methods using exact information.
\end{abstract}

\noindent
{\bf Context:}{\footnotesize
~The material of this report is part of a forthcoming book of the authors on
the evaluation complexity of optimization methods for nonconvex problems.
}

\numsection{Inexact Algorithms Using Dynamic Accuracy}

Most of the literature on optimization assumes that evaluations of the
objective function, as well as evaluations of its derivatives of relevant order(s),
can be carried out exactly.  Unfortunately, this assumption is not always
fulfilled in practice and there are many applications where either the
objective-function values or those of its derivatives (or both) are only known
approximately. This can happen in several contexts. The first is when the
values in questions are computed by some kind of experimental process whose
accuracy can possibly be tuned (with the understanding that more accurate
values maybe be, sometimes substantially, more expensive in terms of
computational effort).  A second related case is when objective-function or
derivatives values result from some (hopefully convergent) iteration:
obtaining more accuracy is also possible by letting the iteration converge
further, but again at the price of possibly significant additional computing.  A
third context, quite popular nowadays in the framework of machine learning, is
when the values of the objective functions and/or its derivatives are obtained
by sampling (say among the terms of a sum involving a very large number of
them). Again, using a larger sample size results in probabilistically better
accuracy, but at a cost.

Extending ideas proposed in \cite{BellGuriMoriToin19}\footnote{For
regularization methods.}, this report discusses a trust-region algorithm which can handle such
contexts, under what we call the ``dynamic accuracy'' requirement: we assume
that \emph{the required values} (objective-function or derivatives) \emph{can
always be computed with an accuracy which is specified, before the
calculation, by the algorithm itself.}  It is also understood in what follows
that \emph{the algorithm should require high accuracy only if necessary},
while \emph{guaranteeing final results to full accuracy}.
In this situation, it is hoped that many function's or derivative's
evaluations can be carried out with a fairly loose accuracy (we will refer to
these as ``inexact values''), thereby resulting in a significantly cheaper
optimization process.

\numsection{Taylor decrements and enforcing accuracy}

We consider the problem of minimizing a smooth, potentially nonconvex,
function $f$ from $\Re^n$ into $\Re$ without constraints on the
variables. This problem has generated a literature too abundant to be reviewed
here, but it is probably fair to say that trust-region methods feature among
the most successful algorithms for its solution, showing excellent practical
performance and solid theoretical background (see \cite{ConnGoulToin00} for an
in-depth discussion).  These methods are based on using n Taylor-series
models, which clearly depend on values and derivatives of the objective
function at a sequence of points (iterates), but in the scenario we are about
to consider, we do not assume that we can calculate them exactly. That is,
rather than having true problem function and derivatives values, $f(x)$ and
$\nabla_x^j f(x)$ for $j \in\ii{q}$ at $x$, we are provided with
approximations $\barf(x)$ and $\overline{\nabla_x^j f}(x)$---here and
hereafter, we denote inexact quantities and approximations with an overbar.

Consequently, while high-degree exact approaches (see
\cite{BirgGardMartSantToin17,BirgGardMartSant20,CartGoulToin20a,BellGuriMoriToin20} for instance) 
deal with a $p$-th degree Taylor-series approximation
\[
T_{f,p}(x,s) = f(x) + \bigsum_{i=1}^p \frac{1}{i!}\nabla_x^i f(x)[s]^i
 \equiv T_{f,p}(x,0) + 
   \bigsum_{i=1}^p \frac{1}{i!} [\nabla_v ^i T_{f,p}(x,v)]_{v=0}[s]^i
\]
of $f$ for perturbations $s$ around $x$, in our new framework, we have to be
content with an inexact equivalent
\[
\barT_{f,p}(x,s) 
 = \barf(x) + \bigsum_{i=1}^p \frac{1}{i!}\overline{\nabla_x^i f}(x)[s]^i
 \equiv \barT_{f,p}(x,0) + 
   \bigsum_{i=1}^p \frac{1}{i!} [\nabla_v ^i \barT_{f,p}(x,v)]_{v=0}[s]^i.
\]
It is therefore pertinent to investigate the 
effect of inexact derivatives on such approximations. As we shall see, 
for our purposes it will be important to achieve sufficient \emph{relative}
accuracy on the value of the Taylor model. More specifically, we will be
concerned with the \emph{Taylor decrement} defined, at $x$ and for a step $s$,
by 
\[
\begin{array}{rl}
\Delta T_{f,p}(x,s) \ngap & \eqdef T_{f,p}(x,0) - T_{f,p}(x,s) \\
& = - \bigsum_{i=1}^p \frac{1}{i!} [\nabla_v ^i T_{f,p}(x,v)]_{v=0}[s]^i \\
& \equiv - \bigsum_{i=1}^p \frac{1}{i!}\nabla_x^i f(x)[s]^i.
\end{array}
\]
While our traditional algorithms depend on this quantity, it is of course of
the question to use them in the present context, as we only have approximate values.
But an obvious alternative is to consider instead the inexact Taylor decrement
\beqn{extDA-barDT}
\barDT_{f,j}(x,s) \eqdef \barT_{f,j}(x,0)-\barT_{f,j}(x,s)
 = - \bigsum_{i=1}^j \frac{1}{i!}
 [\nabla_v^i \barT_{f,j}(x,v)]_{v=0}[s]^i.
 \eeqn
We shall suppose in what follows that a
relative accuracy parameter $\omega \in (0,1)$ is given, and we
will then require that 
\beqn{extDA-barDT-acc}
|\barDT_{f,p}(x,s)-\Delta T_{f,p}(x,s)| \leq \omega \barDT_{f,p}(x,s)
\eeqn
whenever $\barDT_{f,p}(x,s) > 0$. It is not obvious at this point how to enforce this
relative error bound, and we now discuss how this can be achieved.

But Taylor models also occur in termination rule for high-order approximate
minimizers. In particular, it has been argued in \cite{CartGoulToin20a} that ``strong''
approximate $q$-th order minimizers satisfy the necessary optimality condition
\beqn{unco-Taylor-q}
\phi_{f,j}^{\delta_j}(x_\epsilon) \leq \epsilon_j \, \frac{\delta_j^j}{j!}
\;\tim{ for all }\; j\in \ii{q}.
\eeqn
for some $\delta_j\in (0,1)$, where
\beqn{unco-phidef}
\phi_{f,j}^\delta(x)
\eqdef f(x)-\min_{\|d\|\leq\delta}T_{f,j}(x,d),
\eeqn
which is the \emph{largest decrease of the $j$-th order Taylor-series model
$T_{f,j}(x,s)$ achievable by a point at distance at most $\delta$ from
$x$}. Note that $\phi_{f,j}^\delta(x)$  is a continuous function of
$x$ and $\delta$ for $f \in C^j$ \cite[Th. 7]{Hoga73}. It is also important to
observe that $\phi_{f,j}^\delta(x)$ is \emph{independent of the value of}
$f(x_k)$, because the zero-th degree terms cancel in \req{unco-phidef}. In
what follows, we will mostly consider $\delta \leq 1$, but this is not necessary.
Thus $\phi_{f,j}^{\delta_j}(x)$ is itself based on a Taylor-series model and
thus is of importance since we plan to use \req{unco-Taylor-q} as a
termination rule for our proposed algorithm.  This reinforces the need to
understand how to enforce the accuracy which is necessary for the algorithm to
finally produce an exact approximate minimizer.

The attentive reader has noticed that solving the global optimization problem
in \req{unco-phidef}, although not involving any evaluation of $f$ or its
derivatives, still remains a daunting task for $j>2$.  In what follows, we
will allow this calculation to be inexact in the sense that
\req{unco-Taylor-q} will be replaced by the condtion that, for some $\varsigma
\in (0,1]$ and some $d$ with
$\|d\|\leq \delta$,
\beqn{approx-phi-cond}
\varsigma \phi_{f,j}^\delta(x)
\leq \Delta T_{f,j}(x,d)
\leq \varsigma \epsilon_j \, \frac{\delta_j^j}{j!}
\eeqn
for $j\in \ii{q}$. This it to say that, if a given fraction $\varsigma$ of
$\phi_{f,j}^\delta(x)$, the globally optimal Taylor decrease at $x$, can be
calculated, \req{unco-Taylor-q} is still be verifiable at the cost of
reducing the required $\epsilon_j$ by the same fraction.

\subsection{Enforcing the relative error on Taylor decrements}
\label{rel-err-s}

For clarity, we shall temporarily neglect the iteration index $k$.
While there may be circumstances in which \req{extDA-barDT-acc} can be enforced
directly, we consider here that the only control the user has on the accuracy
of $\barDT_{f,j}(x,s)$ is by imposing bounds %$\{\zeta_i\}_{i=1}^j$
on the \emph{absolute} errors of the derivative 
tensors $\{\nabla_x^if(x)\}_{i=1}^j$.
In other words, we seek to ensure \req{extDA-barDT-acc} by selecting
absolute accuracies $\{\zeta_i\}_{i=1}^j$ such that
the desired accuracy requirement follows whenever
\beqn{extDA-vareps-j}
\|\overline{\nabla_x^if}(x)-\nabla_x^if(x)\| \leq \zeta_i
\tim{for} i \in \ii{j},
\eeqn
where $\|\cdot\|$ denotes the Euclidean norm for vectors and the induced
operator norm for matrices and tensors. As one may anticipate by examining
\req{extDA-barDT-acc}, a suitable relative accuracy requirement can be
achieved so long as $\barDT_{f,j}(x,s)$ remains safely away from zero.
However, if exact computations are to be avoided, we may have to 
accept a simpler absolute accuracy guarantee when $\barDT_{f,j}(x,s)$ is 
small, but one that still guarantees our final optimality conditions. 

Of course, not all derivatives need to be inexact in our framework.  If
derivatives of order $i \in \calE\subseteq \ii{q}$ are exact, then the
left-hand side of \req{extDA-vareps-j} vanishes for $i\in \calE$ and
the choice $\zeta_i=0$ for $i\in \calE$ is perfectly adequate. However, we
avoid carrying this distinction in the arguments that follow for the sake of
notational simplicity.

We now start by describing a crucial tool that we use to achieve
\req{extDA-barDT-acc}, the \al{VERIFY} algorithm, inspired by
\cite{BellGuriMoriToin19} and stated as
Algorithm~\ref{extDA-verify} \vpageref[below]{extDA-verify}.
We use this to assess the relative model-accuracy whenever needed in the
algorithms we describe later in this section.

To put our exposition in a general context, we suppose that we 
have a Taylor series $T_r(x,v)$ of a given
function about $x$ in the direction $v$, along with an 
approximation $\barT_r(x,v)$, both of degree $r$, as well as
the decrement $\barDT_r(x,v)$.  
We suppose that a bound $\delta \geq \|v\|$ is given, and that
\emph{required} relative and absolute accuracies $\omega$ and $\xi>0$ 
are on hand. Moreover, we assume that the \emph{current}
upper bounds $\{\zeta_j\}_{j=1}^r$ on absolute accuracies of the derivatives of
$\overline{T}_r(x,v)$ with respect to $v$ at $v=0$ are provided.
Because it will always be the case when we need it, 
we will assume for simplicity that $\barDT_r(x,v) \geq 0$.  Moreover, the
relative accuracy constant $\omega \in (0,1)$ will fixed throughout the forthcoming
algorithms, and we assume that it is given when needed in \al{VERIFY}.

\algo{extDA-verify}{The \al{VERIFY} algorithm}{
\[
\accuracy
= \mbox{\al{VERIFY}}\Big(\delta,\barDT_r(x,v),\{\zeta_i\}_{i=1}^{r},\xi\Big).
\]
\begin{description}
\item[\hspace*{6.5mm}] 
  If
  \vspace*{-3mm}
  \beqn{extDA-verif-term-2}
  \barDT_r(x,v) > 0
  \tim{ and }
  \sum_{i=1}^r \zeta_i \frac{\delta^i}{i!} \leq \omega \barDT_r(x,v),
  \vspace*{-3mm}
  \eeqn
  set \accuracy\ to \relative.
\item[\hspace*{6.5mm}] 
  Otherwise, if
  \vspace*{-3mm}
  \beqn{extDA-verif-term-3}
  \sum_{i=1}^r \zeta_i\frac{\delta^i}{i!} \leq \omega \xi \frac{\delta^r}{r!},
  \vspace*{-3mm}
  \eeqn
  set \accuracy\ to \absolute.
\item[\hspace*{6.5mm}] 
 Otherwise  set \accuracy\ to \insufficient.
  \end{description}
}

\noindent
It will be convenient to say informally that 
\accuracy\ is \sufficient, if it is either \absolute\ or \relative.

We may formalise the accuracy guarantees that result from applying the
\al{VERIFY} algorithm as follows.

\llem{extDA-verify-l}{
Let $\omega \in (0,1]$ and $\delta, \xi$ and $\{\zeta_i\}_{i=1}^{r}>0$. 
Suppose that $\barDT_r(x,v) \geq 0$, that
\[\accuracy
= \mbox{\al{VERIFY}}\Big(\delta,\barDT_r(x,v),\{\zeta_i\}_{i=1}^{r},\xi\Big).
\]
and that
\beqn{extDA-rvareps-j}
\Big\|\Big[\nabla_v^i \barT_r(x,v)\Big]_{v=0}
-\Big [\nabla_v ^i T_r(x,v)\Big]_{v=0}\Big\|
\leq \zeta_i \tim{for} i \in \ii{r}.
\eeqn
Then

\noindent
{\em (i)} \accuracy\ is \sufficient\ whenever
  \vspace*{-3mm}
  \beqn{extDA-max-abs-acc}
    \sum_{i=1}^r\zeta_i\frac{\delta^i}{i!} \leq \omega \xi \frac{\delta^r}{r!},
  \eeqn
 
\noindent
{\em (ii)} if \accuracy\ is \absolute, 
  \beqn{extDA-verif-prop-1}
  \max\Big[\barDT_r(x,v),
    \left|\barDT_r(x,w)- \Delta T_r(x,w)\right|\Big]
  \leq \xi \frac{\delta^r}{r!}
  \eeqn
  for all $w$ with $\|w\|\leq \delta$.

\noindent
{\em (iii)} if \accuracy\ is \relative, 
  $\barDT_r(x,v)>0$ and
  \beqn{extDA-verif-prop-2}
  \left|\barDT_r(x,w)- \Delta T_r(x,w)\right|
  \leq \omega \barDT_r(x,v),
  \tim{for all $w$ with $\|w\|\leq \delta$,} 
  \eeqn

}

\proof{
We first prove proposition \textup{(i)}, and assume that 
\req{extDA-max-abs-acc} holds, which clearly ensures that \req{extDA-verif-term-3} is satisfied. 
Thus either \req{extDA-verif-term-2} or \req{extDA-verif-term-3} must hold and
termination occurs, proving the first proposition.

It follows by definition of the Taylor series, the triangle inequality 
and \req{extDA-rvareps-j} that
\beqn{extDA-difft}
\begin{array}{rl}
\left|\barDT_r(x,w)- \Delta T_r(x,w)\right| \ngap
& = \left| \bigsum_{i=1}^r \bigfrac{
(\nabla^i_w \barT_r(x,w)- \nabla^i_w T_r(x,w))[w]^i}{i!}\right| \\
& \leq \bigsum_{i=1}^r \bigfrac{
\| \nabla^i_w \barT_r(x,w)- \nabla^i_w T_r(x,w)\|\|w\|^i}{i!} \\
& \leq \bigsum_{i=1}^r \zeta_i \frac{\|w\|^i}{i!}.
\end{array}
\eeqn
Consider now the possible \sufficient\ termination cases for the 
algorithm and suppose first that termination occurs with 
\accuracy\ as \absolute.  Then, using \req{extDA-difft},
\req{extDA-verif-term-3} and $\omega < 1$, we have that, 
for any $w$ with $\|w\|\leq \delta$,
\vspace*{-2mm}
\beqn{extDA-DTbarDT-bound}
\left|\barDT_r(x,w)- \Delta T_r(x,w)\right|
\leq \sum_{i=1}^r \zeta_i \frac{\delta^i}{i!}
\leq \omega \xi \frac{\delta^r}{r!}
\leq \xi \frac{\delta^r}{r!}.
\vspace*{-2mm}
\eeqn
If $\barDT_r(x,v) = 0$, we may combine this with \req{extDA-DTbarDT-bound} 
to derive \req{extDA-verif-prop-1}. By contrast, if $\barDT_r(x,v) > 0$,
then since \req{extDA-verif-term-2} failed but \req{extDA-verif-term-3}
holds,
\[
\omega \barDT_r(x,w)
< \sum_{i=1}^r \zeta_i \frac{\delta^i}{i!}
\leq \omega \xi \frac{\delta^r}{r!}.
\]
Combining this inequality with \req{extDA-DTbarDT-bound} yields
\req{extDA-verif-prop-1}. Suppose now that \accuracy\ is \relative.
Then \req{extDA-verif-term-2} holds, and combining it with \req{extDA-difft}
gives that
\vspace*{-2mm}
\[
\left|\barDT_r(x,w)- \Delta T_r(x,w)\right|
\leq \sum_{i=1}^r \zeta_i \frac{\delta^i}{i!}
\leq \omega \barDT_r(x,v_\omega),
\vspace*{-2mm}
\]
for any $w$ with $\|w\|\leq \delta$, which is \req{extDA-verif-prop-2}. 
} %epr

\noindent
Clearly, the outcome corresponding to our initial aim to obtain a relative
error at most $\omega$ corresponds to the case where \accuracy\ is \relative. 
As we will shortly discover, the two other cases are also needed.

\subsection{Computing the approximate optimality measures}\label{computing-barphi-s}

\noindent
Our next concern is how one might compute an optimality measure, given an
inexactly computed $\barDT_{f,p}(x,s)$. Using the crucial measure of optimality 
\beqn{extDA-phiismax}
\phi_{f,j}^\delta(x) \eqdef \max_{\|d\|\leq\delta}\Delta T_{f,j}(x,d),
\eeqn
is out of the question, but an obvious alternative is to consider
instead the {\em inexact measure} 
\beqn{extDA-barphi-def}
\barphi_{f,j}^\delta(x) \eqdef \max_{\|d\|\leq\delta}\barDT_{f,j}(x,d)
\eeqn
that depends on an equivalent sufficiently accurate inexact Taylor decrement.
We immediately observe that $\barphi_{f,j}^\delta(x)$ is independent of the
value of $\barf(x_k)$. Alas, except if we allow exact global minimization in
\req{extDA-barphi-def}, $\barphi_{f,j}^\delta(x)$ may also be to hard to
compute, and we therefore settle to using $\barDT_{f,j}(x,d)$ where $d$ is
such that $\|d\| \leq \delta$ and 
\beqn{whphi-def}
\varsigma \barphi_{f,j}^\delta(x)
\leq \barDT_{f,j}(x,d)
\leq \barphi_{f,j}^\delta(x).
\eeqn

Natural questions are then how well a particular $\barDT_{f,j}(x,d)$ approximates
$\phi_{f,j}^\delta(x)$ and, if there is reasonable agreement, what is a
sensible alternative to the stopping rule \req{unco-Taylor-q}?

We answer both questions in
Algorithm~\ref{extDA-step1} below, which shows one way to compute
$\barDT_{f,j}(x,d)$.  For analysis purposes, this algorithm involves a
counter $i_\zeta$ of the number of times accuracy on the derivatives has been
improved.

\algo{extDA-step1}{Computing $\barDT_{f,j}(x,d)$}{ 
The iterate $x_k$, the index $j\in\ii{q}$ and the radius $\delta_k\in (0,1]$
are given, as well as the constants $\gamma_\zeta \in (0,1)$ and $\varsigma\in (0,1]$.
The counter $i_\zeta$, the relative accuracy $\omega \in (0,1]$ 
and the absolute accuracies bounds $\{\zeta_{i,i_\zeta}\}_{i=1}^q$ are also given. 
\\
\vspace*{-3mm}
\begin{description}
\item[Step 1.1: ] If they are not yet available, compute
     $\{\overline{\nabla_x^if}(x_k)\}_{i=1}^j$ satisfying
     \[
     \|\overline{\nabla_x^if}(x_k)-\nabla_x^if(x_k)\| \leq \zeta_{i,i_\zeta}
        \tim{for} i \in \ii{j}.
     \]
     \vspace*{-4mm}
   \item[Step 1.2: ] Find a displacement $d_{k,j}$ with $\|d_{k,j}\|\leq \delta_k$ such that
     \vspace*{-2mm}
     \beqn{trqeda-barphi-approx}
     \varsigma \barphi_{f,j}^{\delta_k}(x_k) \leq \barDT_{f,j}(x_k,d_{k,j})
%     d_{k,j} = 
%     {\rm arg}\max_{\|d\|\leq \delta_k}\barDT_{f,j}(x_k,d)
    \vspace*{-1mm}
    \eeqn
%     and the corresponding Taylor decrement
%     $\barDT_{f,j}(x_k,d_{k,j})$.
     Compute
     \vspace*{-3mm}
%     \[
    \beqn{trqda-verifyo}
     \mbox{\accuracyj} =
     \mbox{\al{VERIFY}}\Big(\delta_k,\barDT_{f,j}(x_k,d_{k,j}),
             \{\zeta_{i,i_\zeta}\}_{i=1}^j,
                 \half \varsigma \epsilon_j\Big).
     \vspace*{-2mm}
     \eeqn
%     \]
     \vspace*{-4mm}
\item[Step 1.3: ]  If \accuracyj is \sufficient, return $\barDT_{f,j}(x_k,d_{k,j})$.
\end{description}

\begin{description}
\item[Step 1.4: ] Otherwise (i.e.\ if \accuracyj is \insufficient), set
     \vspace*{-2mm}
     \beqn{extDA-decrease-vareps}
     \zeta_{i,i_\zeta+1} = \gamma_\zeta\zeta_{i,i_\zeta}
     \tim{ for } i \in \ii{j},
     \vspace*{-2mm}
     \eeqn
     increment $i_\zeta$ by one and return to Step~1.1.
\end{description}
}

\noindent
Observe that known values of derivatives for $i < j$ may be reused in Step~1.1 if
required. 

We now establish that Algorithm~\ref{extDA-step1} produces
values of the required optimality measures that are adequate in the sense that
either an approximate minimizer is detected or a suitable approximation of the
exact optimality measure is obtained.

\llem{extDA-acc-S1-l}{
If Algorithm~\ref{extDA-step1} terminates within Step~1.3 when 
\accuracyj is \absolute, then
\beqn{extDA-term-q-k}
\phi_{f,j}^{\delta_k}(x_k) \leq \epsilon_j \frac{\delta_k^j}{j!}.
\eeqn
Otherwise, if it terminates with \accuracyj as \relative, then
\beqn{extDA-good-phi}
(1-\omega) \barDT_{f,j}(x_k,d_{k,j})
\leq \phi_{f,j}^{\delta_k}(x_k)
\leq (1+\omega) \barDT_{f,j}(x_k,d_{k,j}).
\eeqn
Moreover, termination with one of these two outcomes must occur if
\beqn{extdA-stop-step1}
\max_{i\in\ii{j}}\zeta_{i,i_\zeta} \leq \frac{\omega}{4}\,\varsigma \epsilon_j\,\frac{\delta_k^{j-1}}{j!}
\eeqn
}%

\proof{
Consider $j\in \ii{q}$. We first notice that Step~1.1 of
Algorithm~\ref{extDA-step1} yields \req{extDA-rvareps-j} with
$T_r=T_r^{f}$ and $r=j$, so that the assumptions of
Lemma~\ref{extDA-verify-l} are satisfied. Note first that, because
of \req{trqeda-barphi-approx} and since
$\barphi_{f,j}^{\delta_k}(x_k)\geq 0$ by definition,
we have that $ \barDT_{f,j}(x_k,d_{k,j})\geq 0$. Suppose now
that the
\al{VERIFY} algorithm returns {\tt accuracy} = {\tt absolute} but now
$\barDT_{f,j}(x_k,d_{k,j})\geq 0$.
%and thus $d_{k,j}\neq 0$.
Using the fact that the nature of Step~1.2 ensures that
$\varsigma \barDT_{f,j}(x_k,d) \leq \barDT_{f,j}(x_k,d_{k,j})$ for $d$ with
$\|d\|\leq\delta_{k,j}$ we have, using \req{extDA-verif-prop-1} with $\xi
= \half \varsigma \epsilon_j$, that, for all such $d$, 
\[
  \begin{array}{lcl}
  \varsigma \Delta T_{f,j}(x_k,d)
  & \leq & \varsigma \barDT_{f,j}(x_k,d) + 
  \varsigma\left|\barDT_{f,j}(x_k,d)- \Delta T_{f,j}(x_k,d)\right| \\*[1ex]
  & \leq & \barDT_{f,j}(x_k,d_{k,j}) + 
  \left|\barDT_{f,j}(x_k,d)- \Delta T_{f,j}(x_k,d)\right| \\*[1ex]
  & \leq & \varsigma \epsilon_j \bigfrac{\delta_{k,j}^j}{j!}
  \end{array}
\]
yielding  \req{extDA-term-q-k}. If the \al{VERIFY} algorithm returns {\tt
accuracy} = {\tt relative},
then, for any $d$ with $\|d\|\leq \delta_{k,j}$,
\[
\begin{array}{lcl}
\Delta T_{f,j}(x_k,d)
& \leq & \barDT_{f,j}(x_k,d) + 
  \left|\barDT_{f,j}(x_k,d)- \Delta T_{f,j}(x_k,d)\right| \\*[1ex]
& \leq & (1+\omega)\barDT_{f,j}(x_k,d_{k,j}).
\end{array}
\]
Thus, for all $d$ with $\|d\|\leq \delta_{k,j}$,
\[
\begin{array}{lcl}
\max\Big[0,\Delta T_{f,j}(x_k,d)\Big]
&  \leq &(1+\omega)\max\Big[0,\barDT_{f,j}(x_k,d_{k,j})\Big]\\*[1ex]
& = & (1+\omega)\barDT_{f,j}(x_k,d_{k,j})
\end{array}
\]
and the rightmost part of \req{extDA-good-phi} follows.
Similarly, for any $d$ with $\|d\|\leq \delta_{k,j}$,
\[
\begin{array}{lcl}
\Delta T_{f,j}(x_k,d)
& \geq & \barDT_{f,j}(x_k,d) - 
\left|\barDT_{f,j}(x_k,d)- \Delta T_{f,j}(x_k,d)\right| \\*[1ex]
& \geq & \barDT_{f,j}(x_k,d) - \omega\barDT_{f,j}(x_k,d_{k,j}).
\end{array}
\]
\vspace*{-2mm}
Hence
\[
\begin{array}{lcl}
\bigmax_{\|d\|\leq \delta_{k,j}}\Delta T_{f,j}(x_k,d)
& \geq & \bigmax_{\|d\|\leq \delta_{k,j}}
\left[ \barDT_{f,j}(x_k,d) -
      \omega\barDT_{f,j}(x_k,d_{k,j})\right] \\*[2ex]
& \geq & (1-\omega)\barDT_{f,j}(x_k,d_{k,j}).
\end{array}
\]
Since $\barDT_{f,j}(x_k,d_{k,j})>0$ when the \al{VERIFY} algorithm returns
{\tt accuracy} = {\tt relative}, we then obtain that, for all $\|d\|\leq \delta_{k,j}$,
\[
\begin{array}{lcl}
\max\Big[0,\max_{\|d\|\leq \delta_{k,j}}\Delta T_{f,j}(x_k,d)\Big]
& \geq & (1-\omega)\barDT_{f,j}(x_k,d_{k,j}), \\
%&  = & (1-\omega)\barphi_{f,j}^{\delta_{k,j}}(x_k),
\end{array}
\]
which is the leftmost part of \req{extDA-good-phi}.
In order to prove the last statement of the lemma, suppose
that \req{extdA-stop-step1} holds.  Then
\[
\sum_{i=1}^j\zeta_{i,i_\zeta}\frac{\delta_k^i}{i!}
\leq \max_{i\in\ii{j}}\zeta_{i,i_\zeta} \sum_{i=1}^j \frac{\delta_k^i}{i!}
\leq (\exp(1)-1) \delta_k \max_{i\in\ii{j}}\zeta_{i,i_\zeta}
\leq 2 \delta_k \max_{i\in\ii{j}}\zeta_{i,i_\zeta}
\leq \half \omega \varsigma \epsilon_j \frac{\delta_k^j}{j!}
\]
and Lemma~\ref{extDA-verify-l} (i) then ensures that the call
to \al{VERIFY} in Step~1.2 returns \accuracyj as \sufficient, causing
Algorithm~\ref{extDA-step1} to terminate in Step~1.3.
}%epr

\noindent
Notice that if we apply Algorithm~\ref{extDA-step1} for all
$j \in \ii{q}$ and each returned \accuracyj is \absolute,
the bound \req{extDA-term-q-k} then ensures that
$x_k$ is an $(\epsilon,\delta_k)$-approximate $q$-th-order minimizer.
If \accuracyj is \relative\ and
\beqn{extDA-bar-term-q}
\barDT_{f,j}(x,d_j) \leq 
 \left(\frac{\epsilon_j}{1+\omega}\right) \frac{\delta_j^j}{j!}
 \tim{ for } j \in\ii{q}.
\eeqn
holds for $x = x_k$ and $\|d_j\|\leq \delta_k$
the same is true because of \req{extDA-good-phi}.  Thus checking
\req{extDA-bar-term-q} is an adequate verification of the $j$-th order
optimality condition. Moreover, the call the \al{VERIFY} in Step~1.1 must
return \relative\ if \req{extDA-bar-term-q} fails.
Importantly, these conclusions do not require that the
$\{\barDT_{f,j}(x_k,d_{k,j})\}_{j=1}^q$ use the same set of 
approximate derivatives for all $j\in\ii{q}$, but merely that their accuracy
is deemed \sufficient\ by the \al{VERIFY} algorithm.

\numsection{The \al{TR$q$DA} algorithm and its complexity}\label{TRqDA-s}

In what follows, we shall first consider a trust-region optimization algorithm,
named \al{TR$q$DA} (the \al{DA} suffix refers to the Dynamic Accuracy
framework) whose purpose is to find a vector $x = x_\epsilon$ for which
\req{extDA-bar-term-q} holds
for some vector of optimality radii $\delta \in (0,1]^q$.  This is 
important as we have just shown (in Lemma~\ref{extDA-acc-S1-l}) that 
any $x_k$ investigated by Algorithm~\ref{extDA-step1} for all
$j \in\ii{q}$ is either directly an ($\epsilon,\delta$)-approximate
$q$-th-order minimizer of $f(x)$ because of \req{extDA-term-q-k} or will be if
\req{extDA-bar-term-q} holds at $x = x_k$ because of \req{extDA-good-phi}.

An initial outline of the \al{TR$q$DA} algorithm is presented 
\vpageref[below]{extDA-TRqDA}.

\algo{extDA-TRqDA}{Trust Region with Dynamic Accuracy \\(TR$q$DA, basic version)}
{
\begin{description}
\item[Step~0: Initialisation.]
  A criticality order $q$, a starting point $x_0$ and an initial trust-region radius
  $\Delta_0$ are given, as well as accuracy levels $\epsilon \in (0,1)^q$ and an initial
  set of bounds on absolute derivative accuracies $\{\zeta_{j,0}\}_{j=1}^q$. The
  constants $\omega$, $\varsigma$, $\vartheta$, $\kappa_\zeta$, $\eta_1$, $\eta_2$, $\gamma_1$, $\gamma_2$, 
  $\gamma_3$ and $\Delta_{\max}$ are also given and satisfy
  \[
  \vartheta \in [\min_{j\in\ii{q}} \epsilon_j, 1],
  \;\;\; \Delta_0\leq \Delta_{\max},\;\;\;
  0 < \eta_1 \leq \eta_2 < 1, \;\;\;
  0< \gamma_1 <\gamma_2< 1 <  \gamma_3,
  \]
  \[
  \varsigma \in (0,1],\;\;\;
  \omega \in \Big(0,\min\big[\half \eta_1, \quarter (1-\eta_2)\big]\Big)
  \tim{and}
  \zeta_{j,0} \leq \kappa_\zeta \tim{for} j \in \ii{p}.
  \vspace*{-1mm}
  \]
  Set $k=0$ and $i_\zeta=0$.
\item[Step~1: Termination test.]
  Set $\delta_k = \min[\Delta_k,\vartheta]$. For $j = 1, \ldots, q $,
  \begin{enumerate}
  \item  Evaluate $\overline{\nabla_x^j f}(x_k)$ and compute
  $\barDT_{f,j}(x_k,d_{k,j})$ using Algorithm~\ref{extDA-step1}.
  \item If
  \beqn{extDA-trq-noterm}
  \barDT_{f,j}(x_k,d_{k,j}) > \left(\frac{\epsilon_j}{1+\omega}\right) \frac{\delta_k^j}{j!},
  \eeqn
  go to Step~2 with $d_{k,j}$, the optimality displacement associated with $\barDT_{f,j}(x_k,d_{k,j})$.
  \end{enumerate}
  If the loop on $j$ finishes, terminate with $x_\epsilon = x_k$
  and $\delta_\epsilon = \delta_k$.
\item[Step~2: Step computation.]
  If $\Delta_k \leq \vartheta$, set $s_k = d_{k,j}$ and $\barDT_{f,j}(x_k,s_k)
  = \barDT_{f,j}(x_k,d_{k,j})$.
  Otherwise, compute a step $s_k$ such that $\|s_k\| \leq \Delta_k$,
  \beqn{extDA-trq-decrease}
  \barDT_{f,j}(x_k,s_k) \geq \barDT_{f,j}(x_k,d_{k,j})
  \eeqn
  and \req{extDA-barDT-acc} holds---see
  Algorithm~\ref{extDA-trqda-step2} below for details.
\item[Step~3: Accept the new iterate.]
   Compute $\barf(x_k+s_k)$ ensuring that
  \vspace*{-2mm}
  \beqn{extDA-trq-Df+-DT}
  |\barf(x_k+s_k)-f(x_k+s_k)| \leq \omega \barDT_{f,j}(x_k,s_k).
  \vspace*{-1mm}
  \eeqn
  Also ensure (by setting
  $\barf(x_k)=\barf(x_{k-1}+s_{k-1})$
  or by recomputing $\barf(x_k)$) that 
  \vspace*{-3mm}
  \beqn{extDA-trq-Df-DT}
  |\barf(x_k)-f(x_k)| \leq \omega \barDT_{f,j}(x_k,s_k).
  \vspace*{-1mm}
  \eeqn
  Then compute
  \beqn{extDA-trq-rhok-def}
  \rho_k = \frac{\barf(x_k) - \barf(x_k+s_k)}
                {\barDT_{f,j}(x_k,s_k)}.
  \eeqn
  If $\rho_k \geq \eta_1$, then set
  $x_{k+1} = x_k + s_k$; otherwise set $x_{k+1} = x_k$.
\item[Step~4: Update the trust-region radius.]
  Set
  \beqn{extDA-trq-Delta-update}
  \Delta_{k+1} \in \left\{ \begin{array}{ll}
  {}[\gamma_1\Delta_k, \gamma_2 \Delta_k] & \tim{if} \rho_k < \eta_1,\\
  {}[\gamma_2\Delta_k, \Delta_k] & \tim{if} \rho_k \in  [\eta_1, \eta_2),\\
  {}[\Delta_k, \min(\Delta_{\max},\gamma_3 \Delta_k)] & \tim{if} \rho_k \geq  \eta_2,
  \end{array}\right.
  \eeqn
  Increment $k$ by one and go to Step~2 with $d_{k+1,j}=d_{k,j}$ if $x_{k+1} = x_k$
  and $\Delta_{k+1}\geq \theta$, or to Step~1 otherwise.
  \end{description}
  }

This algorithm  does not specify how to
find the step required in Step 2. This vital ingredient will be the subject of
what will follow. In addition, we stress that although \req{extDA-trq-Df+-DT}
and \req{extDA-trq-Df-DT} might suggest that we need to know the true $f$,
this is not the case, rather we simply need some mechanism to ensure that 
$x_k$ and $x_k+s_k$ satisfy the required bounds. These bounds are 
needed to guarantee convergence. Notice that the value of
$\barDT_{f,j}(x_k,d_{k,j})$ and  $\barDT_{f,j}(x_k,s_k)$ do not depend on
the value of $\barf(x_k)$, and so Step~1 and 2 are also independent of this
value. In particular, this allows to postpone the choice of $\barf(x_0)$ to
Step~3. At iteration $k$, a new value of $\barf(x_k)$ has to be computed to
ensure \req{extDA-trq-Df-DT} in Step~3 only when
$\barDT_{f,j}(x_{k-1},s_{k-1})>\barDT_{f,j}(x_k,s_k)$. If this is the case,
the (inexact) function value is computed twice rather than once in that
iteration. Finally note that the choice $\vartheta = 1$ is acceptable since we
have assumed that $\epsilon_j \leq 1$ for all $j\in\ii{q}$.

As usual for trust-region algorithms, iteration $k$ is said to be successful
when $\rho_k \geq \eta_1$ and $x_{k+1}= x_k+s_k$, and we define $\calS$,
$\calS_k$ and $\calU_k$ as 
\beqn{unc1-SU-def}
\calS \eqdef \{ k \in \mathbf{N} \mid \rho_k \geq \eta_1 \}
\tim{ and }
\calU \eqdef  \mathbf{N} \setminus \calS,
\eeqn
the sets of \emph{successful} and \emph{unsuccessful} iterations, respectively,
and
\beqn{unc1-SUk-def}
\calS_k \eqdef \{ j \in \iiz{k} \mid \rho_j \geq \eta_1 \}
\tim{ and }
\calU_k \eqdef \iiz{k} \setminus \calS_k,
\eeqn
the corresponding sets up to iteration $k$. Notice that $x_{k+1} = x_k+s_k$
for $k \in \calS$, while $x_{k+1}=x_k$ for $k\in \calU$.

For future reference,
we now state a property of the \al{TR$q$DA} algorithm that solely depends on the
mechanism \req{extDA-trq-Delta-update} to update the trust-region radius.

\llem{SvsU}{
Suppose that the \al{TR1} algorithm is used and that $\Delta_k \geq
\Delta_{\min}$ for some $\Delta_{\min} \in (0,\Delta_0]$. Then
\beqn{unc1-tr1-unsucc-neg}
k \leq |\calS_k| \left(1+\frac{\log\gamma_3}{|\log\gamma_2|}\right)
+ \frac{1}{|\log\gamma_2|}
\left|\log\left(\frac{\Delta_{\min}}{\Delta_0}\right)\right|,
\eeqn
}

\proof{
Observe that \req{extDA-trq-Delta-update} and our assumption imply that
\[
\Delta_{i+1}\leq \gamma_3 \Delta_i, \quad i\in \calS_k
\tim{ and }
\Delta_{i+1}\leq \gamma_2 \Delta_i, \quad i\in \calU_k.
\]
Using our assumption, we thus deduce inductively that
\[
\Delta_{\min}
\leq \Delta_k
\leq \Delta_0\gamma_3^{|\calS_k|}\gamma_2^{|\calU_k|}.
\]
which gives that
\[
\gamma_3^{|\calS_k|}\gamma_2^{|\calU_k|}
\geq \frac{\Delta_{\min} }{\Delta_0}
\]
and we obtain inequality \req{unc1-tr1-unsucc-neg} by taking logarithms on both
sides and recalling that $\gamma_2\in (0,1)$ and that $k = |\calS_k| +
|\calU_k|$.
} % epr

\noindent
In words, so long as the trust-region radius is bounded from below, the
total number of iterations performed thus far is bounded in terms of the
number of successful ones. Note that this lemma is independent of the specific
choice of $s_k$.

\subsection{Computing the step $s_k$}

We now have to specify how to compute the step $s_k$ required by Step~2
whenever $\Delta_k > \vartheta$, in which case $\delta_k = \vartheta$. While any step satisfying both $\|s_k\|\leq
\Delta_k$ and \req{extDA-trq-decrease} is acceptable, we still have to provide
a mechanism that ensures \req{extDA-barDT-acc}.  This is the aim of
Algorithm~\ref{extDA-trqda-step2}.

\algo{extDA-trqda-step2}{Detailed Step 2 of the \al{TR$q$DA} algorithm when $\Delta_k>\vartheta$}{
The iterate $x_k$, the relative accuracy $\omega$, the requested accuracy
$\epsilon_j\in (0,1]^q$, the constants $\gamma_\zeta \in (0,1)$, the counter
$i_\zeta$ and the absolute accuracies $\{\zeta_{j,i_\zeta}\}_{j=1}^q$ are given.
The index $j\in\ii{q}$, the optimality displacement $d_{k,j}$
and the constant $\vartheta \in (0,1]$ are also given such that, by \req{extDA-trq-noterm},
\beqn{extDA-trqda-initial-barDT}
\barDT_{f,j}(x_k,d_{k,j}) > \frac{\epsilon_j}{1+\omega} \frac{\vartheta^j}{j!}.
\eeqn
\begin{description}
\item[Step 2.1: ] If they are not yet available, compute
    $\{\overline{\nabla_x^if}(x_k)\}_{i=1}^j$ satisfying
    \[
    \|\overline{\nabla_x^if}(x_k)-\nabla_x^if(x_k)\| \leq \zeta_{i,i_\zeta}
       \tim{for} i \in \ii{j}.
    \]
    \vspace*{-4mm}
\item[Step 2.2: Step computation.]
    Compute a step $s_k$ such that $\|s_k\|\leq \Delta_k$ and yielding a
    decrease $\barDT_{f,j}(x_k,s_k)$ satisfying \req{extDA-trq-decrease}.
    Compute
    \vspace*{-2mm}
    \beqn{extDA-trqda-VERIFYsmall}
    \begin{array}{l} \hspace*{-3mm}
    \mbox{\accuracys} = \\
    \hspace*{5mm}
    \mbox{\al{VERIFY}}\Big( \|s_k\|,\barDT_{f,j}(x_k,s_k),\{\zeta_{i,i_\zeta}\}_{i=1}^j,
    \bigfrac{\epsilon_j}{4(1+\omega)}\,\Big(\bigfrac{\vartheta}{\max\big[\vartheta,\|s_k\|\big]}\Big)^j\Big).
    \end{array}
    \vspace*{-1mm}
    \eeqn
  \item [Step 2.3:] If \accuracys\ is \relative,
    go to Step~3 of Algorithm~\ref{extDA-TRqDA}
  with the step $s_k$  and the associated $\barDT_{f,j}(x_k,s_k)$. 
\item[Step 2.4: ] Otherwise, set
  \beqn{extDA-trq-decrease-vareps-2}
  \zeta_{i,i_\epsilon+1} = \gamma_\zeta\zeta_{i,i_\epsilon}
  \tim{for} i \in \ii{j},
  \eeqn
  increment $i_\zeta$ by one and go to Step~2.1.
\end{description}
}

The next lemma reassuringly shows that 
Algorithm~\ref{extDA-trqda-step2} must terminate, and provides
useful details of the outcome.

\llem{extDA-trqda-good-step2}{
  Suppose that the detailed Step~2 given by
  Algorithm~\ref{extDA-trqda-step2} is used in the \al{TR$q$DA}
  algorithm whenever $\Delta_k > \vartheta$. If this condition holds,
  the outcome of the call to \al{VERIFY} in Step~2.2 is
  \relative\ and termination must occur with this outcome if
  \beqn{extDA-trqda-stop-step2}
  \max_{i\in\ii{j}}\zeta_{i,i_\zeta}
  \leq \frac{\omega\vartheta^{j-1}}{8j!(1+\omega)} \,\epsilon_j.
  \eeqn.
  In all cases, we have that $\barDT_{f,j}(x_k,s_k) > 0$  and
  \beqn{extDA-relacc-sk-ok}
  \left|\barDT_{f,j}(x_k,s_k)-\Delta T_{f,j}(x_k,s_k)\right|
  \leq \omega \barDT_{f,j}(x_k,s_k).
  \eeqn
  }%

\proof{Suppose first that $\Delta_k \leq \vartheta$.  Then $s_k=d_{k,j}$ and
  \req{extDA-trq-noterm} gives $\barDT_{f,j}(x_k,s_k) > 0$.  Moreover, our
  comment at the end of Section~\ref{computing-barphi-s} shows that the outcome of the
  \al{VERIFY} algorithm called in Step~1.1 for order $j$ must be \relative.
  Lemma~\ref{extDA-verify-l}(iii) then ensures that \req{extDA-relacc-sk-ok}
  holds.
  
  Suppose now that $\Delta_k > \vartheta$ and thus $\delta_k=\vartheta$. We therefore
  have that Algorithm~\ref{extDA-trqda-step2} was used to compute $s_k$.
  Because derivatives may be re-evaluated within the course of this
  algori\-thm, we need to identify the
  particular inexact Taylor series we are considering: we will therefore 
  distinguish $\overline{T}_{f,j}^0(x_k,d_{k,j})$, $\barDT_{f,j}^{0}(x_k,d_{k,j})$ and
  the corresponding accuracy bounds $\{\zeta_i^0\}_{i=1}^j$ using
  the derivatives $\{\overline{\nabla_x^if}(x_k)\}_{i=1}^j$ available on
  entry of the algorithm, from $\overline{T}_{f,j}^+(x_k,d_{k,j})$,
  $\barDT_{f,j}^+(x_k,d_{k,j})$ and $\{\zeta_i^+\}_{i=1}^j$ using derivatives 
  after one or more executions of Step 2.4. By construction, we have that
  \beqn{extDA-trqda-zetas}
  \zeta_i^+ < \zeta_i^0
  \tim{ for } i\in\ii{j}.
  \eeqn
  We also note that, by \req{extDA-trq-noterm} and \req{extDA-trq-decrease},
  \beqn{extDA-trqda-big-barDT0}
  \barDT_{f,j}^0(x_k,s_k^0)
  \geq \barDT_{f,j}^0(x_k,d_{k,j})
  > \frac{\epsilon_j}{1+\omega} \,\frac{\vartheta^j}{j!}>0,
  \eeqn
  where $s_k^0$ is computed using $\overline{T}_{f,j}^0$.

  Observe now that the \al{TR$q$DA} has not terminated at Step~1 and thus 
  that \req{extDA-trq-noterm} holds.  
  This in turn implies that 
  \[
  2 \barDT_{f,j}^0(x_k,d_{k,j})
  > (1 + \omega_k )  \barDT_{f,j}^0(x_k,d_{k,j})
  > \epsilon_j \frac{\vartheta^j}{j!}
  \]
  since $\omega < 1$,   and hence the call the
  \al{VERIFY} in Step~1.2 of Algorithm~\ref{extDA-step1}
  has returned \accuracyj\ as \relative.
  Therefore  \req{extDA-verif-term-2} must hold with $\zeta_i =\zeta_i^0$, $\delta=\vartheta$ and
  $\xi = \half \epsilon_j$, yielding that
  \beqn{extDA-trqda-sumbound}
  \bigsum_{i=1}^j\zeta_j^0\frac{\vartheta^i}{i!}
  \leq \omega \barDT_{f,j}^0(x_k,d_{k,j}).
  \eeqn
  As a consequence, we find that
  \beqn{extDA-trqda-DTlower}
  \begin{array}{lcl}
  \Delta T_{f,j}(x_k,d_{k,j})
  & \geq & \barDT_{f,j}^0(x_k,d_{k,j})
          - |\barDT_{f,j}^0(x_k,d_{k,j})-\Delta T_{f,j}(x_k,d_{k,j})| \\*[2ex]
  & \geq & \barDT_{f,j}^0(x_k,d_{k,j})
          - \bigsum_{i=1}^j\zeta_j^0\frac{\vartheta^i}{i!}\\*[2ex]
  & \geq & (1-\omega)\barDT_{f,j}(x_k,d_{k,j})\\*[2ex]
  &   >  & \bigfrac{1-\omega}{1+\omega}\,\epsilon_j\,\bigfrac{\vartheta^j}{j!}
  \end{array}
  \eeqn
  from the triangle
  inequality, \req{extDA-difft}, the definition of $\{\zeta_i^0\}_{i=1}^j$,
  the fact that $\|d_{k,j}\| \leq \delta_k=\vartheta$ and \req{extDA-trqda-big-barDT0}.
  Using similar reasoning, but now with \req{extDA-trqda-zetas}, 
  we also deduce that 
  \beqn{extDA-trqda-ineq1}
  \begin{array}{lcl}
  \barDT_{f,j}^+(x_k,d_{k,j})
  & \geq & \Delta T_{f,j}(x_k,d_{k,j})
          - |\barDT_{f,j}^+(x_k,d_{k,j})-\Delta T_{f,j}(x_k,d_{k,j})| \\*[2ex]
  & \geq & \Delta T_{f,j}(x_k,d_{k,j})
          - \bigsum_{i=1}^j\zeta_j^+\frac{\vartheta^i}{i!}\\*[2ex]
  & > & \Delta T_{f,j}(x_k,d_{k,j}) -
          \bigsum_{i=1}^j\zeta_j^0\frac{\vartheta^i}{i!}.
  \end{array}
  \eeqn
  Combining this with \req{extDA-trqda-sumbound} and \req{extDA-trqda-DTlower}
  \[
  \begin{array}{lcl}
  \barDT_{f,j}^+(x_k,d_{k,j})
  & \geq & \Delta T_{f,j}(x_k,d_{k,j}) - \omega \barDT_{f,j}^0(x_k,d_{k,j}) \\*[2ex]
  & \geq & \Delta T_{f,j}(x_k,d_{k,j}) - \left(\bigfrac{\omega}{1-\omega} \right)
           \Delta T_{f,j}(x_k,d_{k,j}) \\*[2ex]
  & \geq & \bigfrac{1-\omega}{1+\omega}\left(1 - \bigfrac{\omega}{1-\omega} \right)
           \,\epsilon_j\,\bigfrac{\vartheta^j}{j!}  \\*[2ex]
  & > & \bigfrac{\epsilon_j}{4(1+\omega)} \, \bigfrac{\vartheta^j}{j!},
  \end{array}
  \]  
  where we have used the fact that $\omega < \quarter(1-\eta_2) < \quarter$ to
  deduce the last inequality. Hence, because of \req{extDA-trq-decrease},
  \beqn{extDA-trqda-big-barDT+}
  \barDT_{f,j}^+(x_k,s_k^+)
  \geq \barDT_{f,j}^+(x_k,d_{k,j})
  > \frac{\epsilon_j}{4(1+\omega)} \,\frac{\vartheta^j}{j!}>0.
  \eeqn
  Suppose now that $\barDT_{f,j}(x_k,s_k)$  is any of $\barDT_{f,j}^0(x_k,s_k)$ or
  $\barDT_{f,j}^+(x_k,s_k)$, and
  that the call to
  \al{VERIFY} in \req{extDA-trqda-VERIFYsmall} returns \absolute.
  Applying Lemma~\ref{extDA-verify-l} (ii), we deduce that
  \[
  \barDT_{f,j}(x_k,s_k)
  \leq \frac{\epsilon_j}{4(1+\omega)}\,\frac{\vartheta^j}{\max\big[\vartheta,\|s_k\|\big]^j}\,\frac{\|s_k\|^j}{j!}
  \leq \frac{\epsilon_j}{4(1+\omega)} \,\frac{\vartheta^j}{j!},
  \]
  which contradicts both \req{extDA-trqda-big-barDT0} and
  \req{extDA-trqda-big-barDT+}. This is thus impossible and the 
  call to \al{VERIFY} in \req{extDA-trqda-VERIFYsmall} must also return either
  \relative\ or \insufficient.
  It also follows from \req{extDA-trqda-big-barDT0} and
  \req{extDA-trqda-big-barDT+} that 
  \beqn{ombardt}
  \omega \barDT_{f,j}(x_k,s_k) >
    \frac{\omega \epsilon_j}{4(1+\omega)} \,\frac{\vartheta^j}{j!} > 0,
  \eeqn
  and thus if Step~2.4 continues to be called, ultimately
  \req{extDA-trq-decrease-vareps-2} will ensure that
  \beqn{extDA-trqda-step2-innerbound}
  \frac{\omega \epsilon_j}{4(1+\omega)} \,\frac{\vartheta^j}{j!}
  \geq \sum_{i=1}^j \zeta_{i,i_\zeta} \frac{\vartheta^i}{i!}.
  \eeqn
  This and \req{ombardt} then imply that eventually \req{extDA-verif-term-2}
  in the call to \al{VERIFY} in \req{extDA-trqda-VERIFYsmall} will hold, and
  hence \accuracys\ is \relative.  Thus the exit test in Step 2.3 will
  ultimately be satisfied, and Algorithm~\ref{extDA-trqda-step2}
  will terminate in a finite number of iterations with 
  $\barDT_{f,j}(x_k,s_k) > 0$, because of \req{extDA-trqda-big-barDT0} and
  \req{extDA-trqda-big-barDT+}, and \accuracys\ as \relative.  We
  may then apply Lemma~\ref{extDA-verify-l} (iii) to obtain
  \req{extDA-relacc-sk-ok}.
  Finally observe that, since $\vartheta \leq 1 $, we have that
  \[
  \sum_{i=1}^j \zeta_{i,i_\zeta} \frac{\vartheta^i}{i!}
  \leq \max_{i\in\ii{j}}\zeta_{i,i_\zeta}\sum_{i=1}^j \frac{\vartheta^i}{i!}
  \leq (\exp(1)-1) \vartheta \max_{i\in\ii{j}}\zeta_{i,i_\zeta}
  \leq 2 \vartheta \max_{i\in\ii{j}}\zeta_{i,i_\zeta}.
  \]
  Combining this with \req{extDA-trqda-stop-step2} and using 
  \req{ombardt}, we deduce that \req{extDA-verif-term-2}
  in the call to \al{VERIFY} in \req{extDA-trqda-VERIFYsmall} will hold,
  \accuracys\ is \relative, and termination of 
  Algorithm~\ref{extDA-trqda-step2} in Step 2.3 will occur.
}  % epr

\noindent
The aim of the mechanism of the second item of Step~2.2 should now be clear:
the choice of the last argument in the call to \al{VERIFY} in 
\req{extDA-trqda-VERIFYsmall} is designed to ensure that the outcome
\absolute\ cannot happen.  
This is achieved by ensuring progressively shorter steps are taken 
unless a large inexact decrement is obtained. Observe that the choice $s_k =
d_{k,j}$ is always possible and guarantees that inordinate accuracy is never
needed.

Observe that the mechanism of Algorithm \ref{extDA-trqda-step2}
allows loose accuracy if the inexact decrease $\barDT_{f,j}(x_k,s_k)$ is
large---the test \req{extDA-verif-term-2} will be satisfied in the call to
\al{VERIFY} in Step~2, and thus \al{VERIFY} ignores its last, absolute
accuracy argument ($\xi$) in this case---even if the trust-region radius is
small, while it demands higher absolute accuracy if a large step results in a
small decrease.

For future reference, we note that the last argument in the call to
\al{VERIFY} in \req{extDA-trqda-VERIFYsmall} satisfies
\beqn{extDA-trqda-lastarg}
\frac{\epsilon_j}{4(1+\omega)}\,\frac{\vartheta^j}{\max\big[\vartheta,\|s_k\|\big]^j}
\geq \frac{\epsilon_j}{4(1+\omega)}\,\frac{\vartheta^j}{\max\big[1,\Delta_{\max}\big]^j}
\eeqn

\subsection{Evaluation complexity for the \al{TR$q$DA} algorithm}

We are now ready to analyse the complexity of the \al{TR$q$DA} algorithm of
\vpageref{extDA-TRqDA}, where Step~2 is implemented as in
Algorithm~\ref{extDA-trqda-step2}. We first state our assumptions.

\begin{description}
  \item[AS.1] The function $f$ from $\Re^n$ to $\Re$ is $p$ times continuously
    differentiable and each of its derivatives $\nabla_x^\ell f(x)$ of order
    $\ell\in\ii{p}$ is Lipschitz continuous, that is, for every $j\in\ii{p}$ there exists a constant
    $L_{f,j} \geq 1$ such that, for all $x,y\in\Re^n$,
    \beqn{Lipschitz-f}
    \|\nabla_x^j f(x) - \nabla_x^j f(y)\| \leq L_{f,j} \|x-y\|,
    \eeqn
\item[AS.2] There is a constant $f_{\rm low}$ such that $f(x) \geq f_{\rm
  low}$ for all $x\in \Re^n$.
\end{description}

\noindent
For simplicity of notation, define
\beqn{unco-Lf-def}
L_f \eqdef \max[1, \max_{j\in\ii{q}} L_{f,j}].
\eeqn
The Lipschitz continuity of the derivatives of $f$  has a crucial consequence.

\llem{tech-Taylor-theorem}{
Suppose that AS.1 holds. Then for all $x,s \in \Re^n$,
\beqn{tech-resf}
|f(x+s) - T_{f,j}(x,s)| \leq \bigfrac{L_{f,j}}{(j+1)!} \, \|s\|^{j+1}.
\eeqn
%\beqn{tech-resder}
%\| \nabla^\ell_x f(x+s) -  \nabla^\ell_s T_{f,j}(x,s) \|
%\leq\bigfrac{L_{f,\ell}}{(j-\ell+1)!} \|s\|^{j-\ell+1}.
%\ms (\ell = 1,\ldots, j).
%\eeqn
}

\proof{See \cite[Lemma~2.1]{CartGoulToin20b} with $\beta=1$.}

\noindent
We start our analysis with a simple observation.

\llem{extDA-trqda-phihat}{
At iteration $k$ before termination of the \al{TRqDA} algorithm, 
define \beqn{extDA-trqda-phihat-def}
\sphi_{f,k} \eqdef \frac{j!\,\barDT_{f,j}(x_k,d_{k,j})}{\delta_k^j},
\eeqn
where $j$ is the index for which 
$\barphi_{f,j}^{\delta_k}(x_k) > \epsilon_j / (1+\omega) \delta_k^j/j!$
in Step~1 of the iteration.  Then
\beqn{extDA-trqda-phihat-lower}
\min_{i \in \iiz{k}} \sphi_{f,i}  \geq
\bigfrac{\epsmin}{1+\omega}
\eeqn
where $\epsmin = \min_{j\in\ii{q}} \epsilon_j$.  Moreover,
\beqn{unco-trqda-DTjk}
\barDT_{f,j}(x_k,s_k)
\geq \sphi_{f,k}\frac{\delta_k^j}{j!}
\eeqn
}

\proof{
Let $k$ be the index of an iteration before termination.  Then the mechanism
of Step~1 ensures the existence of $j$ such that
\beqn{unco-trqda-phib}
\barDT_{f,j}(x_k,d_{k,j})
> \frac{\epsilon_j}{1+\omega} \,\frac{\delta_k^j}{j!}
\geq \frac{\epsmin}{1+\omega} \,\left(\frac{\delta_k^j}{j!}\right).
\eeqn
The definiton of $\sphi_{f,k}$ then directly implies that
\beqn{unco-trq-sphi-bound1}
\sphi_{f,k} \geq \frac{\epsmin}{1+\omega}.
\eeqn
Since termination has not yet occurred at iteration $k$, the same inequality
must hold for all iterations $i\in\iiz{k}$, yielding
\req{extDA-trqda-phihat-lower}.  The bound \req{unco-trqda-DTjk}
directly results from
\[
\barDT_{f,j}(x_k,s_k)
\geq \barDT_{f,j}(x_k,d_{k,j})
= \sphi_{f,k}\frac{\delta_k^j}{j!},
\]
where we have used \req{extDA-trq-decrease} to derive the first inequality and
the definitions of $\barphi_{f,j}^{\delta_k}(x_k)$ and $\sphi_{f,k}$ to
obtain the equalities.
} % epr

\noindent
We now derive an ``inexact'' variant of the condition 
that ensures that an iteration is very successful.

\llem{extDA-trqda-condsucc-l}{
Suppose that AS.1 holds, and that
$\sphi_{f,k}$ is defined by \req{extDA-trqda-phihat-def}.
Suppose also that
\beqn{extDA-trqda-condsucc}
\Delta_k
\leq \min\left\{\vartheta,\frac{1-\eta_2}{4\max[1,L_f]}\,\sphi_{f,k} \right\}
\eeqn
at iteration $k$ of Algorithm~\ref{extDA-TRqDA}.
Then $\rho_k \geq \eta_2$,  iteration $k$ is very successful and
$\Delta_{k+1}\geq \Delta_k$.    
}

\proof{
We first note that \req{extDA-trqda-condsucc} implies that $\delta_k =
\min[\vartheta,\Delta_k] = \Delta_k$. Then we may use
\req{extDA-trq-rhok-def}, the triangle inequality, \req{extDA-trq-Df+-DT} and
\req{extDA-trq-Df-DT} and \req{extDA-relacc-sk-ok} (see
Lemma~\ref{extDA-trqda-good-step2}) successively to deduce that 
\[
\begin{array}{ll}
|\rho_k - 1|
& \leq \bigfrac{|\barf(x_k+s_k) - \barT_{f,j}(x_k,s_k)|}
        {\barDT_{f,j}(x_k,s_k)} \\*[2ex]
& \leq \bigfrac{1}{\barDT_{f,j}(x_k,s_k)}
   \Big[|\barf(x_k+s_k) - f(x_k+s_k)| \\*[2ex]
&\hspace*{13.6mm}  + |f(x_k+s_k)-T_{f,j}(x_k,s_k)|
                 + | \barT_{f,j}(x_k,s_k)-T_{f,j}(x_k,s_k)| \Big] \\*[2ex]
& \leq \bigfrac{1}{ \barDT_{f,j}(x_k,s_k)}
   \Big[ |f(x_k+s_k)-T_{f,j}(x_k,s_k)|+3\omega \barDT_{f,j}(x_k,s_k)\Big].\\*[2ex].
\end{array}
\]
Invoking \req{tech-resf} in
Lemma~\ref{tech-Taylor-theorem}, the bound $\|s_k\|\leq \Delta_k=\delta_k$,
\req{unco-Lf-def}, \req{unco-trqda-DTjk}, the fact that $\omega \leq \quarter (1-\eta_2)$, 
and \req{extDA-trqda-condsucc}, we deduce that
\[
|\rho_k - 1|
 \leq \bigfrac{L_{f,j}\,\delta_k^{j+1}}{(j+1)\,\delta_k^j \sphi_{f,k}} + 3\omega
 \leq \bigfrac{L_f\Delta_k}{\sphi_{f,k}} +\sfrac{3}{4}(1-\eta_2)
 \leq 1-\eta_2
 \]
and thus that $\rho_k \geq \eta_2$. Then iteration $k$ is very successful and
\req{extDA-trq-Delta-update} then yields that $\Delta_{k+1} \geq
\Delta_k$. 
}

\noindent
This allows us to derive lower bounds on the trust-region radius
and the model decrease.

\llem{extDA-trq-lemma}{
Suppose that AS.1 holds. 
Then, for all $k\geq0$,
\beqn{extDA-trq-Delta-lower-new}
\Delta_k
\geq  \min\left\{\gamma_1\vartheta,\kappa_{\Delta} \min_{i\in\iiz{k}}\sphi_{f,i}\right\}
\eeqn
where $\sphi_{f,i}$ is defined in \req{extDA-trqda-phihat-def} and 
\beqn{extDA-trq-kappaDdef}
\kappa_{\Delta} \eqdef \frac{\gamma_1(1-\eta_2)}{\max[1,L_f]}
\min\left[\vartheta,\frac{\Delta_0\min_{j\in\ii{q}}\delta_{0,j}^j}
  {2q(\max_{i\in\ii{q}}\|\nabla_x^if(x_0)\|+\kappa_\zeta)}\right]
\eeqn
}

\proof{
Note that, using \req{extDA-trqda-phihat-def}, \req{unco-phidef}, \req{whphi-def}  and the bounds
$\|\overline{\nabla_x^if}(x_0)\| \leq \|\nabla_x^if(x_0)\|+\kappa_\zeta$ and $\delta_{0,j}\leq 1$, we have that
\[
\begin{array}{lcl}
\sphi_{f,0}
& \leq & \bigmax_{j\in\ii{q}}\bigfrac{j!\,\barDT_{f,j}(x_0,d_0)}{\delta_{0,j}^j}\\
%& \leq & \bigmax_{j\in\ii{q}}\bigfrac{j!\,\barphi_{f,j}^{\delta_{0,j}}(x_0)}{\delta_{0,j}^j}\\
& \leq & q \bigmax_{j\in\ii{q}} \bigfrac{\max_{i\in\ii{j}}\|\overline{\nabla_x^if}(x_0)\|}{\delta_{0,j}^j}
   \,\bigsum_{i=1}^j\bigfrac{\delta_{0,i}^i}{i!}\\
& \leq & q\,
   \bigfrac{\max_{i\in\ii{q}}\|\overline{\nabla_x^if}(x_0)\|}{\min_{j\in\ii{q}}\delta_{0,j}^j}\,\big({\rm
     exp}(\delta_{0,j})-1\big)\\
& \leq & 2q\,
   \bigfrac{\max_{i\in\ii{q}}\|\nabla_x^if(x_0)\|+\kappa_\zeta}{\min_{j\in\ii{q}}\delta_{0,j}^j}
\end{array}
\]
and thus, since $\gamma_1(1-\eta_2)<1$,
\[
\kappa_{\Delta}
\leq \frac{\gamma_1(1-\eta_2)}{\max[1,L_f]}
\min\left[\vartheta,\frac{\Delta_0}{\sphi_{f,0}}\right]
\leq \frac{\Delta_0}{\sphi_{f,0}}.
\]
As a consequence, \req{extDA-trq-Delta-lower-new} holds for $k=0$.
Suppose now that $k \geq 1$ is the first iteration such that
\req{extDA-trq-Delta-lower-new} is violated. The updating rule
\req{extDA-trq-Delta-update} then ensures that 
\beqn{unco-trq-new2}
\Delta_{k-1} < \frac{1-\eta_2}{L_f}\,\min_{i\in\iiz{k}}\sphi_{f,i}
\tim{and}
\Delta_{k-1}\leq \vartheta.
\eeqn
Moreover, since           
\[
\sphi_{f,k-1}
\geq \min_{i\in\iiz{k-1}}\sphi_{f,i}
\geq \min_{i\in\iiz{k}}\sphi_{f,i}
\]
we deduce that
\beqn{unco-trq-new3}
\bigfrac{\Delta_{k-1}}{\sphi_{f,k-1}}
\leq \bigfrac{\Delta_{k-1}}{\bigmin_{i\in\iiz{k-1}}\sphi_{f,i}}
\leq \bigfrac{\Delta_{k-1}}{\bigmin_{i\in\iiz{k}}\sphi_{f,i}}
< \bigfrac{1-\eta_2}{L_f}.
\eeqn
Lemma~\ref{extDA-trqda-condsucc-l} and the second part of \req{unco-trq-new2}
then ensure that $\Delta_{k-1}\leq
\Delta_k$.  Using this bound, the second inequality of \req{unco-trq-new3} and
the fact that \req{extDA-trq-Delta-lower-new} is violated at iteration $k$, we
obtain that
\[
\bigfrac{\Delta_{k-1}}
     {\bigmin_{i\in\iiz{k-1}}\sphi_{f,i}}
\leq  \bigfrac{\Delta_k}
     {\bigmin_{i\in\iiz{k}}\sphi_{f,i}}
     < \gamma_1 \bigfrac{1-\eta_2}{L_f}
     \tim{ and }
\Delta_{k-1} \leq \gamma_1 \vartheta.
\]
As a consequence \req{extDA-trq-Delta-lower-new} is also violated at iteration $k-1$.
But this contradicts the assumption that iteration $k$ is the first such that
this inequality fails.  This latter assumption is thus impossible,
and no such iteration can exist. 
}

\llem{extDA-trq-model-decrease}{
Suppose that AS.1 holds. Then, for all $k\geq 0$ before termination,
\beqn{extDA-trq-model-decrease}
\barDT_{f,j}(x_k,s_k)
\geq \frac{\kappa_\delta^{q+1}}{q!} \, \epsmin^{q+1},
\eeqn
where
\beqn{extDA-trqda-kappddef}
\kappa_\delta \eqdef \frac{\kappa_\Delta}{1+\omega}.
\eeqn
}%

\proof{
Suppose first that $\Delta_k > \vartheta$ and therefore $\delta_k=\vartheta$. Then
\req{unco-trqda-DTjk} and the bound $\vartheta \geq \epsmin$ give that
\[
\barDT_{f,j}(x_k,s_k)
\geq \frac{\vartheta^j}{j!}\,\sphi_{f,k}
\geq \frac{\vartheta^q}{q!}\,\sphi_{f,k}
\geq \frac{\epsmin^q}{q!}\frac{\epsmin}{1+\omega},
\]
which yields \req{extDA-trq-model-decrease} since $\kappa_\Delta <1$ and
$\varsigma\leq 1$.
If $\Delta_k \leq \vartheta$, then $\delta_k= \Delta_k$,
and \req{extDA-trq-Delta-lower-new} implies that
\[
\delta_k
\geq \min\left\{\gamma_1\vartheta,\kappa_{\Delta} \min_{i\in\iiz{k}}\sphi_{f,i}\right\}.
\]
Therefore, using \req{unco-trqda-DTjk} again,
\beqn{unco-trq-DT-lower}
\begin{array}{lcl}
\barDT_{f,j}(x_k,s_k)
& \geq & \bigfrac{1}{j!}\sphi_{f,k}\min\left\{\gamma_1\vartheta,\kappa_\Delta \bigmin_{i\in\iiz{k}}\sphi_{f,k}\right\}^j\\
& \geq & \bigfrac{1}{q!}\min\left\{\gamma_1\vartheta,\kappa_\Delta \bigmin_{i\in\iiz{k}}\sphi_{f,k}\right\}^{q+1}.
\end{array}
\eeqn
Moreover, if termination hasn't occurred at iteration $k$,
we have that \req{extDA-trqda-phihat-lower} holds,
and, because $\kappa_\delta \leq 1$ and $\vartheta \geq \epsmin$,
\req{unco-trq-DT-lower} in turn implies \req{extDA-trq-model-decrease}. 
} % epr

\noindent
We may now state the complexity bound for the \al{TR$q$DA} algorithm.

\lthm{extDA-trqda-complexity}
{
Suppose that AS.1 and AS.2 hold. Then there exist positive 
constants $\kata{TRqDA}$, $\katb{TRqDA}$,
$\katc{TRqDA}$, $\katd{TRqDA}$, $\kate{TRqDA}$ and $\katf{TRqDA}$
such that, for any $\epsilon \in (0,1]^q$, the \al{TR$q$DA} algorithm requires
at most 
\beqn{extDA-trqda-fcomp}
\begin{array}{l}
\evbndf{TRqDA}  \eqdef
\kata{TRqDA}\bigfrac{f(x_0)-f_{\rm low}}{\bigmin_{j\in\ii{q}}\epsilon_j^{q+1}}
   +\katb{TRqDA}\left|\log\left(\bigmin_{j\in\ii{q}}\epsilon_j\right)\right|
   +\katc{TRqDA}\\*[5ex]
\hspace*{9mm} = \eo{\bigmax_{j\in\ii{q}}\epsilon_j^{-(q+1)}}
\end{array}
\eeqn
(inexact) evaluations of $f$ and at most
\beqn{extDA-trqda-dcomp}
\begin{array}{l}
\evbndd{TRqDA}  \eqdef
\katd{TRqDA}\bigfrac{f(x_0)-f_{\rm low}}{\bigmin_{j\in\ii{q}}\epsilon_j^{q+1}}
   +\kate{TRqDA}\left|\log\left(\bigmin_{j\in\ii{q}}\epsilon_j\right)\right|
   +\katf{TRqDA}\\*[5ex]
\hspace*{9mm} = \eo{\bigmax_{j\in\ii{q}}\epsilon_j^{-(q+1)}}
\end{array}
\eeqn
(inexact) evaluations of $\{\nabla_x^f\}_{j=1}^q$
to produce an iterate
$x_\epsilon$ and an optimality radius $\delta_\epsilon\in (0,1]$ such that
$\phi_{f,q}^{\delta_\epsilon}(x_\epsilon) \leq \epsilon_j \delta_\epsilon^j/j!$
for all $j \in \ii{q}$,
}

\proof{
If $i$ is the index of a successful iteration before termination, we have that
\beqn{extDA-trqda-decreaseb}
\begin{array}{lcl}
f(x_i)-f(x_{i+1})
& \geq & [\barf(x_i) - \barf(x_{i+1})]
          - 2\omega \barDT_{f,j}(x_i,s_i)\\*[1.5ex]
& \geq & \eta_1 \barDT_{f,j}(x_i,s_i)-  2\omega  \barDT_{f,j}(x_i,s_i)\\*[1.5ex]
& \geq & \bigfrac{(\eta_1-2\omega) \kappa_\delta^{q+1}}{q!}\, \epsmin^{q+1} > 0
\end{array}
\eeqn
using successively \req{extDA-trq-Df+-DT} and \req{extDA-trq-Df-DT}
\req{extDA-trq-rhok-def}, \req{extDA-trq-model-decrease}
and the requirement that $\omega < \half \eta_1$.
Now let $k$ be the index of an arbitrary iteration before termination.
Using AS.2, the nature of successful iterations and
\req{extDA-trqda-decreaseb}, we deduce that 
\[
f(x_0) - f_{\rm low}
\geq  f(x_0) - f(x_{k+1})
= \bigsum_{i\in \calS_k} [f(x_i)-f(x_{i+1}) ]
\geq |\calS_k|\,[\kats{TRqDA}]^{-1}\,\epsmin^{q+1},
\]
where
\beqn{extDA-trqda-kqptrqdas}
\kats{TRqDA} = \frac{q!}{(\eta_1-2\omega)\kappa_\delta^{q+1}},
\eeqn
and thus that the total number of \emph{successful iterations} before
termination is given by
\beqn{extDA-trqda-nsucc}
|\calS_k| \leq \kats{TRqDA}\frac{f(x_0) - f_{\rm low}}
          {\epsmin^{q+1}}.
\eeqn
Observe now that combining respectively \req{extDA-trq-Delta-lower-new}, 
\req{extDA-trq-kappaDdef} and \req{extDA-trqda-phihat-lower},  we obtain that
\beqn{extDA-trqda-Deltabelow}
\Delta_k \geq \min(\vartheta, \Delta_k)
\geq 
\min\left(\vartheta,\kappa_{\Delta} \min_{i\in\iiz{k}}\sphi_{f,i}\right)
= \kappa_{\Delta} \min_{i\in\iiz{k}}\sphi_{f,i}
\geq \kappa_{\delta} \epsmin,
\eeqn
We may then invoke Lemma~\ref{SvsU} to deduce that the total
number of iterations required is bounded by
\[
|\calS_k| \left(1+\frac{\log\gamma_3}{|\log\gamma_2|}\right)
+ \frac{1}{|\log\gamma_2|}
%\left(|\log(\epsmin)| + \left|\log\left(\frac{\kappa_{\delta\!\min}}{\Delta_0}
\left(|\log(\epsmin)| + \left|\log\left(\frac{\kappa_{\delta}}{\Delta_0}
\right)\right|\right)+1.
\]
and hence the total number of approximate function evaluations is at most
twice this number, which yields \req{extDA-trqda-fcomp} with the coefficients
\req{extDA-trqda-kqptrqdas},  
\beqn{extDA-trqda-katab}
\kata{TRqDA} \eqdef 2\kats{TRqDA}\left(1+\frac{\log\gamma_3}{|\log\gamma_2|}\right)
\ms
\katb{TRqDA} \eqdef \bigfrac{2}{|\log\gamma_2|}
\eeqn
and
\beqn{extDA-trqda-katc}
\katc{TRqDA}
\eqdef  \bigfrac{2}{|\log\gamma_2|}\left|\log\left(
\bigfrac{\kappa_{\delta}}{\Delta_0}\right)\right|
 + \bigfrac{2}{|\log(\gamma_\zeta)|} +2.
\eeqn
In order to derive an upper bound on the the number of derivatives'
evaluations,
we now have to count the number of additional derivative
evaluations caused by the need to approximate them to the desired
accuracy. Observe that repeated evaluations at a given iterate $x_k$ are only
needed when the current values  of the absolute errors are smaller than used
previously at $x_k$. These 
absolute errors are, by construction, linearly decreasing with rate $\gamma_\zeta$,
Indeed, they are initialised in Step~0 of the \al{TR$q$DA} algorithm,
decreased each time by a factor $\gamma_\zeta$ in \req{extDA-decrease-vareps}
invoked in Step~1.4 of  Algorithm~\ref{extDA-step1}, down to values
$\{\zeta_{j,i_\zeta}\}_{j=1}^q$ which are then passed to the modified Step~2,
and possibly decreased there further in \req{extDA-trq-decrease-vareps-2}
in Step~2.3 of Algorithm~\ref{extDA-trqda-step2} again by successive
multiplication by $\gamma_\zeta$.
We now use \req{extdA-stop-step1} in Lemma~\ref{extDA-acc-S1-l} and
\req{extDA-trqda-stop-step2} in Lemma~\ref{extDA-trqda-good-step2} to
deduce that the maximal absolute accuracy, $\max_{i\in\ii{j}}\zeta_{i,i_\zeta}$,
will not be reduced below the value
\beqn{minval}
\min\left[ \frac{\omega}{4}\,\varsigma \epsilon_j\,\frac{\delta_k^{j-1}}{j!},
           \frac{\omega}{8(1+\omega)} \,\epsilon_j\,\frac{\delta_k^{j-1}}{j!}\right]
\geq  \frac{\varsigma \omega}{8(1+\omega)} \,\epsilon_j\,\frac{\delta_k^{j-1}}{j!}
\eeqn
at iteration $k$.
But we may now deduce from \req{extDA-trqda-Deltabelow} that
\beqn{extDA-trqda-deltabelow}
\delta_k = \min(\vartheta, \Delta_k)\geq \kappa_{\delta} \epsmin.
\eeqn
This and \req{minval} in turn implies that no further reduction of the 
$\{\zeta_j\}_{j=1}^q$, and hence
no further approximation of $\{\overline{\nabla_x^j f}(x_k)\}_{j=1}^q$, 
can possibly occur in any iteration once the largest initial absolute error
$\max_{j\in\ii{q}}\zeta_{j,0}$ has been reduced by successive multiplications by
$\gamma_\zeta$ sufficiently to ensure that
\beqn{zetamax}
\gamma_\zeta^{i_\zeta} [\max_{j\in\ii{q}}\zeta_{j,0}]
\leq \frac{\varsigma\omega\kappa_{\delta}^{q-1}\epsmin^q}{8(1+\omega)} 
\leq \frac{\varsigma\omega}{8(1+\omega)} \,\epsilon_j\,\frac{(\kappa_{\delta} \epsmin)^{j-1}}{j!}
\leq \frac{\varsigma\omega}{8(1+\omega)} \,\epsilon_j\,\frac{\delta_k^{j-1}}{j!}.
\eeqn
Since the $\zeta_{j,0}$ are initialised in the \al{TR$q$DA} algorithm so that
$\max_{j\in\ii{q}}\zeta_{j,0}\leq \kappa_\zeta$, the bound \req{zetamax}
is achieved once $i_\zeta$, the number of decreases in $\{\zeta_j\}_{j=1}^q$, is
large enough to guarantee that
\[
\gamma_\zeta^{i_\zeta} \kappa_\zeta
\leq  \kap{acc} \epsmin^q
\tim{where} \kap{acc} \eqdef 
 \frac{\varsigma\omega\kappa_{\delta}^{q-1}}{8(1+\omega)},
\]
which is equivalent to asking
\beqn{extDA-trqda-ieps-bound}
i_\zeta \log(\gamma_\zeta) + \log\left(\kappa_\zeta\right)
\leq q\log\left(\epsmin\right) + \log(\kap{acc}).
\eeqn
We now recall that Step~1 of the \al{TR$q$DA} algorithm is only used (and
derivatives evaluated) after successful iterations.
As a consequence, we deduce that the number of evaluations of the derivatives 
of the objective function that occur during the course of the \al{TR$p$DA}
algorithm before termination is at most
\beqn{evald-trqda}
|\calS_k| + i_{\zeta,\min},
\eeqn
i.e., the number iterations in \req{extDA-trqda-nsucc} plus
\[
\begin{array}{lcl}
i_{\zeta,\min}
& \!\! \eqdef \!\! & \left\lfloor
\bigfrac{1}{\log(\gamma_\zeta)}
\left[q\log\left(\epsmin \right)+\log\left(\bigfrac{\kap{acc}}{\kappa_\zeta}\right)\right]
\right\rfloor \\*[2ex]
& \!\!\leq \!\! &
\bigfrac{q}{|\log(\gamma_\zeta)|}\left|\log\left(\epsmin\right)\right|
  + \bigfrac{1}{|\log(\gamma_\zeta)|}
    \left|\log\left(\bigfrac{\kap{acc}}{\kappa_\zeta}\right)\right|+1,
\end{array}
\]
the smallest value of $i_\zeta$ that ensures \req{extDA-trqda-ieps-bound}.
Adding one for the final evaluation at termination, this leads to the desired
evaluation bound \req{extDA-trqda-dcomp} with the coefficients
\[
\katd{TRqDA} = \kats{TRsDA},
\ms
\kate{TRqDA} \eqdef \frac{q}{|\log\gamma_\zeta|}
\tim{and}
\katf{TRqDA} \eqdef
 \bigfrac{1}{|\log(\gamma_\zeta)|}
    \left|\log\left(\bigfrac{\kap{acc}}{\kappa_\zeta}\right)\right|+2.
\]
} % epr

%\absent{%%%%%%%%%%%%%%%%%%%%%%%%%%%%%%%%%%%%%%%%%%%%%%%%%%%%%%%%%%%%%%%%%%%%%%%%

\numsection{Discussion of the \al{TR$q$DA} algorithm}

In order to further avoid overloading notation and over-complicating 
arguments, we have
made a few simplifying assumptions in the description of the \al{TR$q$DA}
algorithm.  The first is that, when accuracy is tightened in Steps~1.4 and 2.4,
we have stipulated a uniform improvement for all derivatives of orders one to
$q$. A more refined version of the algorithm is obviously possible in
which the need to improve accuracy for each derivative is considered 
separately, and that
requires sufficient accuracy on each of the approximate derivatives
$\{\overline{\nabla_x^j f}(x_k)\}_{j=1}^q$. Assuming that $\|s_k\| \leq 1$
and remembering that $\|d_{k,j}\|\leq \delta_k \leq \vartheta \leq 1$, we might
instead consider imposing derivative-specific absolute accuracy requirements
\beqn{extDA-jwise-acc-1}
\|\overline{\nabla_x^\ell f}(x_k) - \nabla_x^\ell f(x_k)\|
\leq \frac{\omega}{3\|s_k\|^\ell}
\barDT_{f,p}(x_k,s_k),
\ms (\ell \in \ii{q}),
\eeqn
and
\beqn{extDA-jwise-acc-2}
\|\overline{\nabla_x^\ell f}(x_k) - \nabla_x^\ell f(x_k)\|
\leq \frac{\omega}{3\|d_{k,j}\|^\ell}
\barDT_{f,j}(x_k,d_{k,j}),
\ms
(j \in \ii{q}, \ell \in \ii{j}),
\eeqn
rather than \req{extDA-barDT-acc} applied to the
directions $s_k$ and $d_{k,j}$  (remember
that they are the only directions used in the \al{VERIFY} tests in
Algorithms~\ref{extDA-step1} and \ref{extDA-trqda-step2}.)
One can then verify that \req{extDA-barDT-acc} still holds for suitable $x$
and $s$ in the computation of $m_k(s_k)$ and $\barDT_{f,j}(x_k,d_{k,j})$.
To see this, consider the accuracy of the Taylor series for $f$ evaluated at
a general step $s$, where $s$ is either $s_k$ or $d_{k,j}$. Using
\req{extDA-jwise-acc-1} or \req{extDA-jwise-acc-2}, we have that, for any $j\in\ii{q}$,
\[
\begin{array}{lcl}
|\barDT_{f,j}(x_k,s)-\Delta T_{f,j}(x_k,s)|
&\leq & \bigsum_{\ell=1}^j \bigfrac{\|s\|^\ell}{\ell!}
        \|\overline{\nabla_x^\ell f}(x_k) - \nabla_x^\ell f(x_k)\|\\*[2ex]
&\leq & \bigsum_{\ell=1}^j \bigfrac{\omega}{3\ell!}
        \barDT_{f,j}(x_k,s)\\*[2ex]
&\leq & \bigfrac{1}{3}\left(\bigsum_{\ell=1}^j\frac{1}{\ell!}\right)\omega\barDT_{f,j}(x_k,s)\\*[2ex]
&\leq & \bigfrac{1}{3}\left(\bigsum_{i=0}^j\frac{1}{i!}\right)\omega\barDT_{f,j}(x_k,s)\\*[2ex]
& <   & \omega \barDT_{f,j}(x_k,s).
\end{array}
\]
Thus \req{extDA-jwise-acc-1} and \req{extDA-jwise-acc-2} guarantee that
\req{extDA-barDT-acc} holds both for $x=x_k$ and $s=s_k$ when $\|s_k\|< 1$
and for $x=x_k$ and $s=d_{k,j}$, as occurring in the computation
of $\barDT_{f,j}(x_k,d_{k,j})$.

Of course, the detailed ``derivative by derivative'' conditions
\req{extDA-jwise-acc-1} and \req{extDA-jwise-acc-2} make no attempt to exploit
possible balancing effects between terms of different degrees $\ell$ in the
Taylor-series model and, in that sense, are more restrictive than
\req{extDA-barDT-acc}. However they illustrate an important point: since the
occurrence of small $\|s_k\| < 1$ and $\|d_{k,j}\|$ can be
expected to happen overwhelmingly often when convergence occurs, the above
conditions indicate that the accuracy requirements on derivatives become
looser for higher-degree derivatives.  This is reminiscent of the situation where a quadratic is
minimized using conjugate-gradients with inexact products, a situation for
which various authors \cite{SimoSzyl03,vandSlei04,GratSimoTitlToin20} have shown
that the accuracy of the products with the Hessian may be progressively relaxed without
affecting convergence.

A second simplifying feature of the \al{TR$q$DA} algorithm relates 
to the insistence  that the absolute accuracies $\{\zeta_{i,i_\zeta}\}_{i=1}^q$ 
are initialised in Step~0 once and for all. As a consequence, the accuracy
requirements can only become more severe as the iteration proceeds.  
This might well be viewed
as inefficient because the true need for accurate derivatives
depends more on their values at a individual rather than the evolving set 
of iterates.  
A version of the algorithm for which
the $\{\zeta_{i,i_\zeta}\}_{i=1}^q$ are reinitialised at every successful
iterate is of course possible, at a moderate increase in the overall
complexity bound.
Indeed, in such a case, the number of ``additional'' derivatives evaluations
$i_{\zeta,\max}$ (in the proof of Theorem~\ref{extDA-trqda-complexity}) 
would no longer need to cover all iterations, but only what happens 
at a single iterate. Thus the logarithmic term in $\epsmin$ is no longer 
added to the number of successful iterations, but
multiplies it, and the worst-case evaluation complexity for the modified
algorithm becomes $\eo{|\log(\epsmin)|\epsmin^{-(q+1)}}$.

\numsection{Conclusions}

We have presented an inexact trust-region algorithm using high-order models
and capable of finding high-order strong approximate minimizers.  We have then
shown that it will find such a $q$-th order approximate minimizer in at most
$\eo{\min_{j\in\ii{q}}\epsilon_j^{-(q+1)}}$ inexact evaluations of the
objective function and its derivatives.  Obviously, the results presented also
cover the case when the function and derivatives evaluations are exact (and
$\omega$ can be set to zero).

{\footnotesize
\section*{\small Acknowledgements}
The authors are grateful to S. Bellavia, G. Gurioli and B. Morini for useful discussions.

%\bibliography{/home/pht/bibs/refs}
%\bibliographystyle{plain}

}

\end{document}